\documentclass[12pt]{article}

\usepackage{microtype} 

\usepackage{graphicx}
\usepackage{amsmath}%
\usepackage{mathtools} 
\usepackage{amsthm}
\usepackage{amsfonts}
\usepackage{amssymb}

\usepackage[cm]{fullpage} 
\usepackage{color}
\usepackage{latexsym}
\usepackage{dsfont}

\usepackage{subfigure}
\usepackage{soul}
\usepackage{enumerate}
\usepackage{multirow}
\usepackage{thm-restate}
\usepackage{xfrac}

\usepackage{diagbox}
\usepackage{makecell}
\usepackage{array,booktabs}
\usepackage{hhline}
\usepackage[table]{xcolor}

\usepackage[numbers,sort,compress,square]{natbib}
\usepackage{doi}

\usepackage{comment}

\newtheorem{theorem}{Theorem}

\newtheorem{corollary}[theorem]{Corollary}
\newtheorem{lemma}[theorem]{Lemma}

\newtheorem{proposition}[theorem]{Proposition}
\theoremstyle{remark}

\theoremstyle{definition}
\newtheorem{example}[theorem]{Example} 
 
\newtheorem{remark}[theorem]{Remark}

\newenvironment{Proof}[1]{\vspace{10pt}\par\noindent{\bf Proof{#1}:}\hspace{3pt}}{}
\newcommand{\qedForProof}{\qed\vspace{10pt}} 

\newcommand{\old}[1]{{{}}}

\graphicspath{{./Figures/}}

\usepackage{algorithm}
\usepackage{algpseudocode}
\def\RHS{right-hand side }

\def\Cov{\textnormal{Cov}}
\def\Exp{\textnormal{Exp}}

\def\M{M}
\def\pGamma{\gamma}

\usepackage{fancyvrb}

\usepackage{float}
\usepackage{comment}

\usepackage{xfrac}
\usepackage{amsmath}%
\usepackage{mathtools} 
\usepackage{amssymb}
\usepackage{makecell}
\usepackage{multirow}

\newcommand\Item[1][]{%
  \ifx\relax#1\relax  \item \else \item[#1] \fi
  \abovedisplayskip=0pt\abovedisplayshortskip=0pt~\vspace*{-\baselineskip}}

\def\Thanks#1{\gdef\thefootnote{\arabic{footnote}}\thanks{#1}}

\usepackage{xurl}

\begin{document}

\title{Power and Limits of Subset Selection in Statistical Estimation}
\author{D. Barak-Pelleg\Thanks{Department of Computer Science, Sami Shamoon College of Engineering, Beer Sheva 8410802, Israel.
E-mail: dina.barak.pelleg@gmail.com}
\and
D.~Berend\Thanks{Department of Mathematics and Institute for the Theory of Computing, Ben-Gurion
University, Beer Sheva 84105, Israel.
E-mail: berend@bgu.ac.il}
}
\maketitle

\begin{abstract}
We study the power and limitations of subset selection in statistical estimation through the framework of \emph{super-teaching}, where a teacher selects a subset of i.i.d.~data to optimize a learner’s estimator. Unlike prior work focused on specific distributions or fixed subset sizes, we develop a general theory under minimal assumptions. 

For mean estimation, we prove that super-teaching is possible for any distribution whose density is bounded away from zero in some neighborhood of the mean, allowing subset sizes growing as $k = o(n^{1/3})$ and achieving error on the order of roughly  $k!/n^{k}$. This significantly extends existing results on admissible distributions and subset scaling. We also extend the analysis to parameters expressed as smooth functionals of expectations, such as variance and scale parameters in classical parametric families, including settings with heavy tails. Moreover, we show that super-teaching can greatly improve estimation rates for nonlinear estimators like the sample median, achieving rates beyond classical asymptotics.

Through examples, including cases where maximum likelihood estimators are inconsistent or fail to be asymptotically normal, we demonstrate that super-teaching can succeed even when standard statistical guarantees break down. Our results establish a unified theory of data selection to enhance statistical efficiency.

\end{abstract}

\section{Introduction}

A central challenge in statistical learning is using data as efficiently as possible. In many modern applications, data sets are large and often contain significant redundancy. Therefore, it is desirable to select a small subset of observations that can preserve or even enhance the statistical performance achieved with the full sample. This approach forms the foundation for various methods, including sample compression, coreset construction, and active learning, all of which aim to identify the most informative examples.

This paper studies a framework introduced by Ma et al.~\cite{ma2018teacher},
known as \emph{super-teaching}. In this setting, we consider a learner who aims to estimate a parameter of an unknown distribution using a prescribed estimation rule. A teacher observes a fixed i.i.d.\ sample from the distribution and may select a subset of the observations to present to the learner. We assume that the teacher has knowledge of the underlying distribution (or of the target parameter). The learner has access only to the selected subset, and not to the full sample. The teacher is not permitted to modify the data, generate synthetic observations, or communicate side information; the only available action is to choose which observations to reveal. The learner then applies the estimation rule to this subset.

The central question is whether subset selection can outperform the full sample in statistical estimation performance.

A key feature of this setting is that the teacher knows the learner’s estimation rule. As a result, the choice of the subset depends on the estimator, and different estimation procedures require different selection strategies.

Super-teaching was introduced by \cite{ma2018teacher} in the setting of mean
estimation under a normal distribution.  
The full sample is of size $n$, but the teacher is allowed to pass only a $k$-subset to the
learner; the learner uses the sample mean.  
The surprising finding is that the teacher can choose a subset whose average
is much closer to the true mean than the full-sample mean.
Specifically, if $k$ is fixed and $n$ is large, one expects among the
$\binom{n}{k}$ subsets at least one whose mean lies roughly between
$(k/n)^{k+\varepsilon}$ and $(k/n)^{k-\varepsilon}$ of the true mean. By contrast, the empirical mean based on the full sample has error of order $n^{-1/2}$.

This result has, however, several important limitations: it is restricted to the normal distribution, considers only mean estimation, and assumes a fixed subset size~$k$. In~\cite{barak2024super}, we extended the result to symmetric unimodal distributions with bounded support and allowed $k$ to grow as $k=o(n^{1/4})$.
Moreover, the error has been reduced to roughly $k!/n^k$.
Although this substantially broadens the scope of the theory, these assumptions still exclude many standard statistical models.

The present paper develops a general theory of super-teaching under minimal assumptions. When estimating the mean, we show that it suffices for the underlying distribution to have a density bounded below in a neighborhood of the mean, without imposing symmetry, unimodality, or bounded support. Under this condition, we obtain super-teaching guarantees for subset sizes growing as $k=o(n^{1/3})$, with error of roughly $k!/n^k$. This significantly enlarges both the class of distributions and the range of subset sizes for which super-teaching is possible.

Beyond mean estimation, our results extend to a broad class of estimators that are smooth functions of empirical means. This includes, for example, variance and scale parameters in many classical parametric families. It also applies in irregular settings, such as heavy-tailed models and cases in which the full-sample maximum likelihood estimator may fail to be consistent.

Section~\ref{sec:main-results} presents the main results, including the general super-teaching theorem and several illustrative examples. Section~\ref{sec:proofs} is devoted to the proofs.

\section{Main Results}
\label{sec:main-results}

Our first main result shows that super-teaching holds under very mild
regularity conditions. Namely, it is sufficient for the density of the underlying distribution to be
bounded below on some  neighborhood of its mean~$\mu$.  
We make no assumptions on symmetry, unimodality, tail behavior, or global
boundedness of the density, as in \cite{barak2024super}.

We use the standard notation $[n]=\{1,2,\ldots,n\}$.

\begin{theorem}
\label{thm:general-result}
Let $X_1,\dots,X_n \stackrel{\mathrm{iid}}{\sim} F$, where $F$ is a continuous
distribution with density~$f$ and mean~$\mu$.  
Assume that $f$ is bounded away from $0$ in a neighborhood of $\mu$.
For any nonempty $I\subseteq[n]$, let
\[
M_I = \frac{1}{|I|}\sum_{i\in I} X_i .
\]
If $k=o(n^{1/3})$, then for every $\delta>0$ and $\varepsilon>0$ there exists
$N$ such that
\begin{equation*}
P\!\left(
\exists\, I\subseteq[n],\ |I|=k:
|M_I - \mu|
\le
\frac{k!}{n^{\,k-\varepsilon}\sqrt{k}}
\right)
\ge 1-\delta ,\qquad n\ge N.
\end{equation*}
\end{theorem}

In other words, under the conditions in the theorem, if the learner estimates the mean by the sample mean, then super-teaching is possible with high probability\footnote{Recall that a sequence of events $(E_j)_{j=1}^\infty$ occurs \emph{with high probability} if  
$\lim_{j\to\infty} P(E_j)=1$.} (w.h.p.). That is, the teacher can provide a set $I\subset [n]$ of size $k$, such that the mean of the observations in the corresponding sub-sample is very close to $\mu$.

The theorem strictly generalizes the results of
\cite{ma2018teacher} and \cite{barak2024super} as indicated above, and in addition enlarges the allowable
growth of $k$ from $o(n^{1/4})$ to $o(n^{1/3})$.  
The bound $k!/n^{\,k}$ is consistent with the combinatorial intuition that the
teacher is effectively choosing the best among $\binom{n}{k}\approx n^k/k!$
options.

\begin{remark}
The restriction $k=o(n^{1/3})$ arises from our second-moment analysis.
We do not know whether it reflects a genuine limitation of the bound in
Theorem~\ref{thm:general-result} or merely a limitation of our proof technique.

At a heuristic level, larger values of $k$ may be expected to facilitate
super-teaching. Indeed, the number of candidate subsets
$\binom{n}{k}$ increases up to $k=\lfloor n/2\rfloor$, while the mean of
a random $k$-subset becomes more concentrated around $\mu$ as $k$
grows. This suggests that the optimal attainable error may continue to
decrease well beyond the range covered by Theorem~\ref{thm:general-result}. However, such
heuristics do not imply that the specific bound of Theorem \ref{thm:general-result} remains
valid for larger $k$. Determining the largest growth rate of $k$ for
which the conclusion of Theorem \ref{thm:general-result} holds remains an interesting open
problem.
\end{remark}

While Theorem \ref{thm:general-result} is phrased in terms of expectation estimation, it may well be applicable to other parameters.

\begin{example}
Consider the classical Cauchy distribution with known scale, say with density
\begin{equation}\label{eq:cauchy-density}
f(x \mid \theta) = \frac{1}{\pi [1 + (x-\theta)^2]}, \qquad x \in \mathbb{R},
\end{equation}
where the location parameter $\theta\in\mathbb R$ is unknown.
Suppose the learner is naive; he estimates $\theta$ by the sample mean
\[
\widehat\theta_{\mathrm{naive}}
=
\frac{1}{n}\sum_{i=1}^n X_i.
\]
Since the Cauchy distribution has no finite mean, this estimator does not
converge at all. (In fact, it is Cauchy distributed with the same parameters as the observations themselves \cite[p.173]{Feller1970}.)

A super-teacher can, however, remedy the situation and push the learner to an excellent
estimate as follows.
Fix any $\alpha>0$ and suppose the teacher first removes all observations outside $(\theta-\alpha,\theta+\alpha)$; for concreteness, take $\alpha=1$.
Since
\[
\mathbb P\bigl(|X_i-\theta|<1\bigr)
=
\frac{1}{2},
\]
for large $n$ the teacher will remain with approximately 
\(n/2\)
observations with high probability.
These observations are drawn from the conditional distribution
$\mathrm{Cauchy}(\theta,1)\mid |X-\theta|<1$, which (by symmetry) has mean~$\theta$ and whose density
is bounded away from 0 throughout $(\theta-1,\theta+1)$. (In fact, the teacher does not really need to discard extreme observations. The theorem guarantees the existence of a subset with exceptional averaging properties within the sub-sample; by using the whole sample the teacher may well find even better sub-samples.)
By Theorem~\ref{thm:general-result}, for any $k=o(n^{1/3})$, with probability
tending to one there exists a sub-sample $(X_i)_{i\in I}$ of size $|I|=k$ such that
\[
\left|
\frac{1}{k}\sum_{i\in I} X_i - \theta
\right|
\le
\frac{k!}{n^{\,k-\varepsilon}\sqrt{k}}.
\]
\end{example}

Many parameters of practical interest can be expressed as smooth functionals of expectations. 
Specifically, suppose the target parameter admits a representation of the form
\[
\theta^* = g\!\left( \mathbb{E}[h(X_i)] \right),
\]
where $h:\mathbb{R}\to\mathbb{R}$ is a fixed measurable function and 
$g:\mathbb{R}\to\mathbb{R}$ is continuously differentiable in a neighborhood of 
$\mathbb{E}[h(X_i)]$. 
In this setting, the natural estimator based on a sub-sample 
\((X_i)_{i \in I}\) is~$ g(\frac{1}{|I|}\sum_{i\in I} h(X_i))$.

We formalize this observation below.

\begin{proposition}
\label{prop:scale-general}
Let  $X_1,\ldots, X_n$  be an i.i.d.\ sample from some distribution and $Y_i = h(X_i)$, $1\le i\le n$. Assume $Y_i$ has a continuous density   bounded away from $0$ in a neighborhood of the mean $m=\mathbb{E}[Y_i]$.  
Let $\theta^*=g(m)$, where $g$ is  continuously differentiable in a neighborhood of~$m$. Denote 
\[
\hat{\theta}_I = g\!\left( \frac{1}{|I|}\sum_{i\in I} Y_i \right), \qquad \emptyset\ne I\subseteq [n].
\]
Then  for every $\delta>0$, $\varepsilon>0$, and $k=o(n^{1/3})$,
for all sufficiently large~$n$,
\begin{equation*}
P\!\left(
\exists\, I\subseteq[n],\ |I|=k:
\left|
\hat{\theta}_I - \theta^*
\right|
\le
\frac{k!}{n^{\,k-\varepsilon}\sqrt{k}}
\right)
\ge 1-\delta.
\end{equation*}
\end{proposition}

As explained after Theorem \ref{thm:general-result}, the implication to super-teaching is straightforward.

\begin{corollary}
Let $X_1,\ldots,X_n$ be an i.i.d.\ sample from a continuous
distribution with mean $\mu$ and variance~$\sigma^2$.
Assume that the density of $X_i$ is bounded away from $0$
in a neighborhood of at least one of the points
$\mu-\sigma$ and $\mu+\sigma$.
For every nonempty $I\subseteq[n]$, define
\[
\hat{\sigma}_I^2
=
\frac1{|I|}
\sum_{i\in I}(X_i-\mu)^2.
\]
Then for every $\delta>0$, $\varepsilon>0$, and
$k=o(n^{1/3})$, for all sufficiently large $n$,
\[
P\!\left(
\exists I\subseteq[n],\ |I|=k:
\left|\hat{\sigma}_I^2-\sigma^2\right|
\le
\frac{k!}{n^{k-\varepsilon}\sqrt{k}}
\right)
\ge 1-\delta .
\]
\end{corollary}

Indeed, apply Proposition~\ref{prop:scale-general} with $h(x)=(x-\mu)^2$ and $g(x)=x$.
Under the stated assumption, the random variable $(X_i-\mu)^2$
has a continuous density bounded away from $0$ in a neighborhood
of its mean~$\sigma^2$.

Proposition~\ref{prop:scale-general} applies to many classical models in which the parameter of interest can be expressed in terms of some moment. In such cases, the teacher can first find a highly accurate subset for the underlying moment and then transfer that accuracy to the target parameter.
\newpage

\begin{example} We illustrate this mechanism for two familiar parametric families:
\begin{enumerate}
\item Let $X_i \sim \mathrm{Gamma}(\alpha,\beta)$ with known shape parameter $\alpha>0$ and unknown scale parameter $\beta>0$. Since $\mathbb{E}[X_i]=\alpha\beta$, the scale parameter can be written in the form $\beta=g(\mathbb{E}[X_i])$, where $g(m)=m/\alpha$. Similarly, when $\beta$ is known, we can write $\alpha = g(\mathbb{E}[X_i])$ with $g(m)=m/\beta$, and obtain a super-teachable estimator of $\alpha$. The same applies if $X_i \sim \mathrm{Inv\text{-}Gamma}(\alpha,\beta)$ (see \cite[p.255]{hoff2009first}) by passing to $Y_i=1/X_i$. 

\item Let $X_i \sim \mathrm{Rayleigh}(\sigma)$ (see \cite[p.169]{papoulis2002probability}) with an unknown scale parameter $\sigma>0$. Then $X_i^2\sim\Exp(1/2\sigma^2)$, so that we can get super-teaching as in the previous case.
\end{enumerate}
\end{example}

Next we consider the problem of estimating the median $M$ of a distribution using the sample median. 
The sample median $\hat{M}$ is a classical example of an L-estimator~\cite[Chap.~21]{vdvaart1998}.
Under mild regularity conditions on the density $f$, in particular continuity and positivity at $M$, standard results on L-estimators imply asymptotic normality for the sample median \cite[Cor.21.5]{vdvaart1998}:
\[
\sqrt{n} \bigl( \hat{M} - M \bigr) \xrightarrow[n\to\infty]{\mathcal{D}} N\left(0, \frac{1/4}{f(M)^2}\right).
\]
Thus, in location families $F(x) = G(x - \theta)$, where $G$ has median $0$, the sample median estimates the location parameter $\theta$ with $\sqrt{n}$-rate of convergence~\cite{vdvaart1998}.

Our next result shows that super-teaching can yield a significant improvement. Given numbers $x_1,\dots,x_n$ and another number $d$, a \emph{most symmetric pair
around $d$} is a pair $x_i,x_j$ such that
\[
\left|\frac{x_i+x_j}{2}-d\right|
=
\min_{1\le i',j'\le n}
\left|\frac{x_{i'}+x_{j'}}{2}-d\right|.
\]
(The indices $i,j$ are allowed to coincide.)

\begin{theorem}
\label{thm:median}
Let $X_1,\dots,X_n \stackrel{\mathrm{iid}}{\sim} F$, where $F$ is a
continuous distribution, and assume its density is bounded away from
$0$ in a neighborhood of its median $M$.
Then for every $\delta>0$ and $\varepsilon>0$ there exists a constant
$N$ such that, for all $n\ge N$,
\[
P\!\left(
\exists\, i,j \in [n] :
\left|
\frac{X_i + X_j}{2} - M
\right|
\le
\frac{1}{n^{2-\varepsilon}}
\right)
\ge 1-\delta .
\]
\end{theorem}

Note that no generality is lost by restricting attention to subsets
of size at most two. Indeed, the median of every sub-sample is either
one of the observations or the midpoint of two observations, and hence
can always be realized as the median of a subset of size one or two.
Thus, based on a sample of size $n$, the teacher can induce the learner
to estimate the median only by one of the $n$ observations or one of
the $\binom{n}{2}$ averages of two of them---altogether $O(n^2)$
possible values. This suggests that, unlike
Theorem~\ref{thm:general-result} and
Proposition~\ref{prop:scale-general}, one should not expect errors much
smaller than order $n^{-2}$.

Theorem~\ref{thm:median} establishes super--teaching with error of
order $O(n^{-2+\varepsilon})$ for every $\varepsilon>0$, which is only
slightly weaker than the natural benchmark $O(n^{-2})$. We do not know
whether the factor $n^{\varepsilon}$ is necessary.

\begin{example}
For the Cauchy location family (scale $1$), the median estimates $\theta$ with asymptotic standard-deviation of $\pi/(2\sqrt{n})$, but super-teaching boosts it to $O(1/n^{2-\varepsilon})$ via the most symmetric pair in some interval around $\theta$, say $[\theta-1,\theta+1]$.
\end{example}

Ma et al \cite[p.1373]{ma2018teacher} conjecture that, if the MLE  satisfies the asymptotic normality condition, it is super-teachable. They also give examples where the condition is not satisfied and super-teaching is impossible. We certainly believe that super-teaching is possible when the MLE satisfies the asymptotic normality condition. Now we provide several examples showing that the opposite is far from true. Namely, super-teaching may be possible even when the condition is not satisfied, and indeed even when the MLE is not consistent.

\begin{example}
We consider the uniform location model \(U(\theta-1/2,\theta+1/2)\) and notice that it behaves quite differently from the scale family $U(0,\theta)$, studied in \cite{ma2018teacher}.
Let \(
X_1,\dots,X_n \stackrel{\mathrm{iid}}{\sim} U(\theta-1/2,\theta+1/2)
\), where $\theta\in\mathbb R$ is unknown. 

A natural learner in this setting uses the estimator
\begin{equation}\label{eq:MLE_uniform}
    \hat\theta_n=\frac{X_{(1)}+X_{(n)}}{2}.
\end{equation}
This estimator is a maximum likelihood estimator, although it is not unique:
every point in the interval
\(\bigl[X_{(n)}-1/2,\;X_{(1)}+1/2\bigr]\)
attains the same likelihood. It is a natural symmetric choice,
depending on the extreme order statistics.

When applied to the full sample, $\hat\theta_n$ is unbiased and consistent, but does not satisfy the asymptotic normality condition. In fact,
\[
2n(\hat\theta_n-\theta)
=
n\bigl(X_{(1)}+X_{(n)}-2\theta\bigr)
\xrightarrow[n\to\infty]{\mathcal{D}} Y_1-Y_2,
\]
where $Y_1$ and $Y_2$ are independent $\mathrm{Exp}(1)$ random variables.
In
particular, after normalization at rate $n$, the limiting distribution is
Laplace rather than Gaussian. (The appearance of a non-Gaussian limit is not
essential here: even if the limit were Gaussian, the estimator would still fail
to be asymptotically normal due to the non-$\sqrt n$ rate.)
Despite this non-regular behavior, if the learner uses the MLE as in~\ref{eq:MLE_uniform}, we have super-teachability. The density is
bounded away from $0$ on the entire support:
\(f(x)=1\) for all $x\in(\theta-\tfrac12,\theta+\tfrac12)$.
Thus the local density condition of Theorem~\ref{thm:general-result} is satisfied.
Although the learner’s estimator is the average of the extreme observations, and not a sample average as in Theorem~\ref{thm:general-result}, when restricted to a
subset of size $k=2$ it coincides with the sample average. Its improvement
under teaching follows then from the same subsample-averaging mechanism captured by the
theorem.
Applying Theorem~\ref{thm:general-result} with $k=2$, the teacher selects a subsample
consisting of two observations $\{x_i,x_j\}$ and provides it to the learner, who uses in turn the average. (Alternatively, we can note that $\theta$ is also the median of the distribution and use Theorem~\ref{thm:median}.) 
Thus, taking $\hat\theta=(X_i + X_j)/2$, the theorem guarantees that, w.h.p., for every $\varepsilon>0$ we have
$|\hat\theta-\theta| \le n^{-2+\varepsilon}$.
\end{example}

\begin{example}
Consider the Cauchy distribution with unknown location parameter $\theta \in \mathbb{R}$ and known scale parameter $1$ (see \eqref{eq:cauchy-density}).

Finding the MLE in this case is technically difficult. The score equation can be written as a polynomial equation of degree $2n-1$. The number of spurious local maxima converges in distribution to $\mathrm{Pois}(1/\pi)$, so with probability about $1-\exp(-1/\pi)\approx 0.273$ the likelihood has multiple local maxima as $n\to\infty$ \cite{reeds1985asymptotic,bai1987maximum}.
These spurious local maxima arise in the tails of the
data and diverge to $\pm\infty$ almost surely as $n$ tends to infinity
\cite{perlman1983limiting}. Consequently, numerical methods such as Newton’s
algorithm are sensitive to the choice of starting values and may converge to
spurious roots; several initialization strategies (e.g., median or truncated
mean) have been proposed and studied
\cite{reeds1985asymptotic,barnett1966evaluation}.

Notwithstanding these technical difficulties, the  MLE does possess
desirable properties in this setting. In particular, recent work shows that the
MLE is consistent and asymptotically normal,
yielding estimation error of order $O(1/\sqrt{n})$
\cite{okamura2026asymptotics}.

Let us explain why $\theta$ is super-teachable. With high probability, there are
$\Omega(n)$ observations in the interval $[\theta-1,\theta+1]$. Suppose first
that the teacher provides the learner with any two such observations, say
$X_i,X_j$. By \cite[Ex.~8.3]{young2005essentials}, the MLE based on these
observations is given by
\[
\hat{\theta} = \frac{X_i + X_j}{2}.
\]
As in the preceding example, by choosing $X_i$ and $X_j$ to be the two
observations that are most symmetric about $\theta$, the learner obtains an
estimate satisfying
\[
|\hat{\theta} - \theta| = O(n^{-2+\varepsilon}).
\]
Thus, even in this heavy-tailed model, super-teaching improves the estimation error from the standard $O(1/\sqrt{n})$ rate for the MLE to the faster rate $O(n^{-2+\varepsilon})$.

\end{example}

\bigskip

Our next example is based on an example of Radford \cite{NealRadford2008MLEBlog}.
\begin{example}
Let
$U\sim U(0,1)$ and $U_{\theta}\sim  U(\theta - e^{-1/\theta^2},\theta)$, where the parameter $ \theta$ belongs to $(0,1)$.
Consider the distribution function
\begin{equation*}
\begin{aligned}
   F(x; \theta)
   &={1}/{2}\cdot F_U(x) + {1}/{2}\cdot F_{U_{\theta}}(x).
\end{aligned}
\end{equation*}
In other words, we toss a coin, and based on the result of this toss draw from $F_U$ or from $F_{U_{\theta}}$.
The corresponding density function is given by:
\begin{equation*}
\begin{aligned}
   f(x; \theta)
   &=\dfrac{1}{2}\cdot\mathds{1}_{[0,1]}(x)
   +\dfrac{e^{1/\theta^2}}{2}\cdot\mathds{1}_{[\theta - e^{-1/\theta^2},\, \theta]}(x)
   =\dfrac{1}{2}
\left(1+{e^{1/\theta^2}}\cdot\mathds{1}_{[\theta - e^{-1/\theta^2},\, \theta]}(x)\right),\qquad 0\le x\le 1.
\end{aligned}
\end{equation*}
Consider the likelihood function for a sample:
\begin{equation}\label{eq:example.neal.1}
\begin{aligned}
  L(\theta;x_1,\ldots, x_n)&= \dfrac{1}{2^n}
\prod_{i=1}^n\left(1+{e^{1/\theta^2}}\cdot\mathds{1}_{[\theta - e^{-1/\theta^2},\, \theta]}(x_i)\right)\\
&\le
\dfrac{1}{2^n}
\left(1+{e^{1/\theta^2}}\right)^{n}\\
&\le
\dfrac{1}{2^n}
\left(2\cdot{e^{1/\theta^2}}\right)^{n}\\
&\le
{e^{n/\theta^2}}.
\end{aligned}
\end{equation}
Denote by $Y'$ the number of observations that fall in $[0, \log n/n]$. Clearly, for sufficiently large $n$,
\(Y'\sim B(n, \log n/2n)\).
By Chernoff's bound, we have at least $\log n/3$ 
observations in $[0,\log n/n]$ w.h.p.
In particular, the minimal observation $x_{\min}$ lies in $[0,\log n/n]$ w.h.p.
Therefore, the  likelihood function at $x_{\min}$ satisfies w.h.p.:
\begin{equation*}
\begin{aligned}
  L(x_{\min};x_1,\ldots, x_n)
&\ge
\dfrac{1}{2^n}
\left(1+{e^{1/x_{\min}^2}}\right)\\
&\ge
\dfrac{e^{n^2/\log^2 n}}{2^n}.
\end{aligned}
\end{equation*}
On the other hand, for ${\theta}'\in[1/{n}^{1/3},1]$, we have
by (\ref{eq:example.neal.1}), 
for large enough $n$: 
\[L({\theta}';x_1,\ldots, x_n)\le e^{n/(1/n^{1/3})^2}=e^{n^{5/3}}<
\dfrac{e^{n^2/\log^2 n}}{2^n}\le L(x_{\min};x_1,\ldots, x_n).\]
It follows that the MLE is at most $1/n^{1/3}$, and in particular tends to $0$ as $n\to\infty$. Thus, the MLE is inconsistent for the distribution in question.

Suppose now that the teacher gives the learner  only one observation, say $x_i$. The likelihood function~is:
\begin{align*}
 L(\theta;x_i)=
 \begin{cases}
     \dfrac{1}{2}\left(1+{e^{1/\theta^2}}\right),
     &\qquad \theta-e^{-1/\theta^2}\le x_i\le \theta,\\[8pt]
     \dfrac{1}{2},&\qquad \text{otherwise}.
 \end{cases}
\end{align*}
(Note that the independent variable is $\theta$; we write it this way as it is impossible to specify the limits for $\theta$ in terms of $x_i$ by elementary functions.)
Since the expression $(1+{e^{1/\theta^2}})$ decreases as a function of $\theta$, the MLE in this case is simply
$\hat\theta = x_i.$
It follows that, for general $n$, by providing the learner with the single observation $x_j$ closest to $\theta$, the teacher usually brings him to an error of $O(1/n)$. (Of course, it may be the case that, by choosing appropriately several observations, the error may become much smaller.)
\end{example}

\section{Proofs}
\label{sec:proofs}

\begin{lemma}
\label{lem:M_I-uniform}
If $X_i\sim U(-1,1)$, $1\le i\le n$, are independent, then for all $k\ge 1$
\begin{equation}
\label{eq:lemma}
P\!\left(
\frac{1}{k}
\left|
\sum_{i=1}^k X_i
\right|
\le \theta
\right)
\ge \sqrt{k}\theta,
\qquad
0<\theta\le \frac{1}{\sqrt{3k}} .
\end{equation}
\end{lemma}

\begin{Proof}{}
The distribution of $\frac{1}{k}\sum_{i=1}^{k}X_i$
is supported on $[-1,1]$ and is symmetric and unimodal 
(see \cite[Sec.D.7.4]{luo2024Bates}
and
\cite[Sec.26.9]{johnson1995continuous}).
It follows that the left-hand side of \eqref{eq:lemma}, namely $p(\theta)=P\left(\left\vert\frac{1}{k}\sum_{i=1}^{k}X_i\right\vert\le \theta\right)$, is concave as a function of~$\theta$ in $[0,1]$. Indeed, if $f$ denotes the density of $\frac1k\sum_{i=1}^k X_i$, then
\[
p(\theta)=2\int_0^\theta f(x)\,dx.
\]
Since $f$ is unimodal and symmetric, $f$ is nonincreasing on
$[0,\infty)$, so $p$ is concave.
On the other hand, the right-hand side of~\eqref{eq:lemma} is linear in~$\theta$ and coincides with   $p(\theta)$ at 0.  
Therefore, it suffices to prove~\eqref{eq:lemma} at the right end point $\theta_0 = 1/\sqrt{3k}$.
By Berry-Esseen's Theorem \cite{Shevtsova2011BerryEsseen}:
\begin{equation*}
\begin{aligned}
1-p(\theta_0)&=1-P\left(\dfrac{1}{k}\left\vert\sum_{i=1}^{k}X_i\right\vert\le \theta_0\right)\\
&
=P\left(\left\vert\frac{\sqrt{3}}{\sqrt{k}}\sum_{i=1}^{k}X_i\right\vert> \sqrt{3k}\theta_0\right)\\
&=
P\left(\left\vert\frac{\sqrt{3}}{\sqrt{k}}\sum_{i=1}^{k}X_i\right\vert> 1\right)
\\
&=2P\left(\frac{\sqrt{3}}{\sqrt{k}}\sum_{i=1}^{k}X_i< -1\right)\\
&=2\left(1-\Phi(1)\right)+\dfrac{2C\rho}{\sigma^3\sqrt{k}},
\end{aligned}
\end{equation*}
where $\Phi$ is the standard normal distribution function, $\sigma=\sqrt{V(X_1)}=1/\sqrt{3}$, $\rho=E(\vert X_1\vert^3)=1/4$, and $C$ is some constant which is less than $0.4748$.
Thus,
\begin{equation*}
\label{eq:lemmaCenter.uniform.2}\begin{aligned}
{p}(\theta_0)&= 2\Phi(1)-1-{2C\rho}/(\sigma^3\sqrt{k})\\
&\ge 2\Phi(1)-1-5/(4\sqrt{k})\\
&=\sqrt{3k}\theta_0\left(0.6826-5/(4\sqrt{k})\right).
\end{aligned}
\end{equation*}
For $k\ge 142$, this gives:
\begin{equation*}
\label{eq:lemmaCenter.uniform.2}\begin{aligned}
{p}(\theta_0)& \ge \sqrt{k}\theta_0 \cdot\sqrt{3}\left(0.6826-5/(4\sqrt{142})\right)\\
& \ge \sqrt{k}\theta_0.
\end{aligned}
\end{equation*}

We verified computationally, using \emph{Python} and \emph{Mathematica}, that~\eqref{eq:lemma} also holds for all
$1 \le k \le 141$. See
\url{https://github.com/dina-barak/power-and-limits-of-subset-selection} for the code and results.

\qedForProof
\end{Proof}

\begin{lemma}
\label{lem:Uniform-Result}
Theorem~\ref{thm:general-result} holds in the special case where $X_1,\dots,X_n \stackrel{\mathrm{iid}}{\sim} U(-1,1)$.
\end{lemma}

The point of proving this very special case of Theorem~\ref{thm:general-result} separately is that we will basically reduce the general case to it.

\begin{Proof}{\ of Lemma~\ref{lem:Uniform-Result}}
Denote by $\theta=\theta(n,k)$ a small number depending on $n$ and $k$.
We want to show that with probability at least $1-\delta$
there exists a subset $I$ of size $k$
such that $|M_I|\le\theta$.

Let $T_I=1$ if $M_I\in\left[-\theta,\theta\right]$, and $T_I=0$ otherwise. Let $T$ be the number of subsets $I$ with $T_I=1$:
\begin{equation*}
    T = \sum_{\vert I\vert = k}T_I.
\end{equation*}
Denote $ \pGamma = P\left(M_I\in\left[-\theta,\theta\right]\right)$.
By Lemma \ref{lem:M_I-uniform},
\begin{equation}\label{eq:new.var.0}
    E(T)
    = \sum_{\vert I\vert = k}E(T_I)
    =\binom{n}{k}\pGamma
    \ge \binom{n}{k}\sqrt{k}\theta.
\end{equation}
Next,
\begin{equation}\label{eq:new.var.1}
\begin{aligned}
V(T)
    &= \sum_{\vert I\vert = k}V(T_I)
    + \sum_{ \vert I\vert = \vert I'\vert = k:\ I\ne I'}\Cov (T_I,T_{I'})\\[0.6em]
    &= \binom{n}{k}V(T_{I_1})
    +\binom{n}{k}\sum_{I':\ 1\le |I_1\cap I'|\le k-1}\Cov(T_{I_1},T_{I'})\\[0.6em]
    &= \binom{n}{k}V(T_{I_1})
    + \binom{n}{k}
    \sum_{j=1}^{k-1}\sum_{I':\  |I_1\cap I'|=j}\Cov(T_{I_1},T_{I'})\\[0.6em]
    &= \binom{n}{k}V(T_{I_1})
    + \binom{n}{k}
    \sum_{j=1}^{k-1}\binom{k}{j}\binom{n-k}{k-j}\Cov(T_{I_1},T_{I'_{\!j}}),
\end{aligned}
\end{equation}
where $I_1=\lbrace1,2,\ldots, k\rbrace$ and $I'_{\!j}=\lbrace 1,2,\ldots,j, k+1,k+2,\ldots, 2k-j\rbrace$.
We start with the variance in the first addend on the \RHS of (\ref{eq:new.var.1})
\begin{equation}\label{eq:new.var.2}
    V(T_{I_1})
= \gamma(1-\gamma) \le \gamma.
\end{equation}
Now we estimate $\Cov(T_{I_1},T_{I'_{\!j}})$:
\begin{align}\label{eq:new.var.2.5}
    \Cov(T_{I_1},T_{I'_{\!j}})
    \le E(T_{I_1}T_{I'_{\!j}})
    =P\left(M_{I_1},M_{I'_{\!j}}\in\left[-\theta,\theta\right]\right).
\end{align}
Denote $J={I_1}\cap I'_{\!j}=\{1,\ldots, j\}$ and recall that
$$M_J = \frac{1}{\vert J\vert}\sum_{i\in J}X_i = \frac{1}{j}\sum_{i=1}^jX_i,$$
Suppose that $M_J$ assumes the value $a$, namely, the sum of the $X_i$-s over $i\in J$ is $a\vert J\vert=a\cdot j$.
The  events $M_{I_1}\in\left[-\theta,\theta\right]$ and $M_{I'_{\!j}}\in\left[-\theta,\theta\right]$  will occur if and only if the sum of the $X_i$-s over all $i\in {I_1}-J=\{j+1,j+2,\ldots, k\}$ falls in 
$$\left[-k\theta -a j ,\ k\theta -a j\right],$$
and the same happens with the sum of the $X_i$-s over all $i\in I'_{\!j}-J=\{k+1,k+2,\ldots, 2k-j\}$. Denote by $f_J$ the density function of the variable $M_J$. 
As the $M_J$-s are supported on $[-1,1]$, and as ${I_1}-J$ and $I'_{\!j}-J$ are disjoint, 
\begin{equation}\label{eq:new.var.3}
\begin{aligned}
   &P\left(M_{I_1},M_{I'_{\!j}}\in\left[-\theta,\theta\right]\right)\\
   &\hspace{10mm}=\int\limits_{-1}^{1}
   f_J(t)\cdot P\left(\left(k-j\right)M_{{I_1}-J}\in\left[-k\theta-t j, k\theta -tj\right]
   \cap\left(k-j\right)M_{I'_{\!j}-J}\in\left[-k\theta-tj, k\theta -tj\right]\right)dt\\
   &\hspace{10mm}=\int\limits_{-1}^{1}
   f_J(t)\cdot P\left(\left(k-j\right)M_{I_1-J}\in\left[-k\theta-tj, k\theta -tj\right]\right)^2dt.
\end{aligned}
\end{equation}
Consider the probability in the second factor of the integrand. The event in question occurs when $M_{I_1-J}$ is in a certain interval of length $2k\theta/(k-j)$. 
As mentioned in the proof of Lemma~\ref{lem:M_I-uniform},
$M_{I_1-J}$ is symmetric and unimodal. Therefore, this probability is maximal when the interval is symmetric with respect to $0$. Hence, among all intervals of a given length, the interval centered at $0$ has maximal probability. Therefore,
\begin{equation}\label{eq:new.var.4}
\begin{aligned}    &P\left(\left(k-j\right)M_{I_1-J}\in\left[-k\theta-tj, k\theta -tj\right]\right)\\
    &\hspace{10mm}\le P\left(M_{I_1-J}\in\left[\frac{-k\theta}{k-j}, \frac{k\theta }{k-j}\right]\right)
    =P\left(\sum_{i=1}^{ k-j}X_i \in\left[{-k\theta}, {k\theta }\right]\right)\\
&\hspace{10mm}=P\left(\sum_{i=1}^{k-j}\dfrac{X_i+1}{2} \in\left[{\dfrac{-k\theta+k-j}{2}}, {\dfrac{k\theta+k-j}{2} }\right]\right).
\end{aligned}
\end{equation}
Note that, for a positive integer $r$, the sum $\sum_{i=1}^{r}(X_i+1)/2$ 
is a sum of $r$ i.i.d.\ $U(0,1)$-distributed variables. (Its distribution is the so-called Irwin-Hall distribution with parameter $r$.) As its density function is bounded above by $1$ for every $r$, we have
\begin{equation}\label{eq:new.var.5}
\begin{aligned}
    &P\left(\left(k-j\right)M_{I_1-J}\in\left[-k\theta, k\theta \right]\right)\le \dfrac{k\theta+k-j}{2}-\dfrac{-k\theta+k-j}{2}
    = k\theta.
\end{aligned}
\end{equation}
By (\ref{eq:new.var.2.5})-(\ref{eq:new.var.5})
\begin{equation}\label{eq:new.var.6}
\begin{aligned}
   \Cov(T_{I_1},T_{I'_{\!j}})
   &\le (k\theta)^2\cdot\int\limits_{-1}^{1}
   f_J(t)dt=(k\theta)^2.
\end{aligned}
\end{equation}
Consider the products of binomial coefficients on the \RHS of (\ref{eq:new.var.1}):
\begin{equation}
\label{eq:new.var.6.0}
\begin{aligned}
 \binom{k}{j}
 \binom{n-k}{k-j}
 &\le\dfrac{k^j}{j!}\cdot\dfrac{(n-k)^{k-j}}{(k-j)!}
 \le\dfrac{k^{2j}\cdot n^{k-j}}{j!\,k!}.
 \end{aligned}
\end{equation}
Thus, by (\ref{eq:new.var.6}) and since $k=o(n^{1/3})$, the second sum  on the \RHS of (\ref{eq:new.var.1})
is:
\begin{equation}\label{eq:new.var.6.1}
\begin{aligned}
 \sum_{j=1}^{k-1}
 \binom{k}{j}
 \binom{n-k}{k-j}\Cov(T_{I_1},T_{I'_{\!j}})
 &\le
\sum_{j=1}^{k-1}\dfrac{k^{2j}}{j!}\cdot\dfrac{n^{k-j}}{k!}\cdot {k^2}\theta^2\\
 &\le
\sum_{j=1}^{k-1}\dfrac{k^{2j+2}}{n^{j-1}}\cdot\dfrac{n^{k-1}}{k!}\cdot \theta^2\cdot \dfrac{1}{j!}\\
&\le
\dfrac{k^{4}\cdot n^{k-1}\cdot \theta^2}{k!}
\sum_{j=1}^{k-1}\dfrac{1}{j!}\\
&\le
ek^{4}\cdot n^{k-1}\cdot \dfrac{\theta^2}{k!}.
\end{aligned}
\end{equation}
Thus, by (\ref{eq:new.var.1}), (\ref{eq:new.var.2}), and (\ref{eq:new.var.6.1}),
\begin{equation}\label{eq:new.var.7.0}
\begin{aligned}
V(T)    &\le\binom{n}{k}\gamma
    + e\binom{n}{k}k^{4}\cdot n^{k-1}\cdot \dfrac{\theta^2}{k!}.
\end{aligned}
\end{equation}
By the second moment method,
\begin{equation*}
\begin{aligned}
    P(T=0)&= 1- P(T>0)
    \le 1-\frac{E^2\left(T\right)}{E\left(T^2\right)}
=\frac{V\left(T\right)}{E\left(T^2\right)}
    \le \frac{V\left(T\right)}{E^2\left(T\right)}.
\end{aligned}
\end{equation*}
Thus, by (\ref{eq:new.var.0}), (\ref{eq:new.var.7.0}) and Lemma~\ref{lem:M_I-uniform}, taking $\theta = k!/n^{k-\varepsilon}\cdot 1/{\sqrt{k}}$, we obtain
\begin{equation}\label{eq:new.var.9.1}
\begin{aligned}
    P(T=0)&
    \le \frac{\binom{n}{k}\gamma
    + e\binom{n}{k}
k^{4}\cdot n^{k-1}
\cdot\theta^2/k!}
    {\left(\binom{n}{k}\gamma\right)^2}\\[0.5em]
&=\frac{1}
    {\binom{n}{k}}
    \left(\frac{1}
    {\gamma}
    +\frac{
ek^{4}\cdot n^{k-1}}
    {\gamma^2}\cdot\dfrac{\theta^2}{k!}\right)  \\[0.5em]
&\le \frac{k!}
    {n(n-1)\cdots(n-k+1)}
    \left(\frac{1}
    {\sqrt{k}\theta}
    +\frac{
ek^{4}\cdot n^{k-1}}
    {k\cdot k!}\right)\\[0.5em]
&= \frac{k!}
    {n^k\cdot(1+o(1))}
    \left(\frac{1}
    {\sqrt{k}\cdot  k!/n^{k-\varepsilon}\cdot1/{\sqrt{k}}}
    +\frac{
ek^{3}\cdot n^{k-1}}
    {k!}\right)\\[0.5em]
    &= \frac{1}{n^{\varepsilon}}\left(1+o(1)\right)
    +O\left({\frac{k^3}{n}}\right).
\end{aligned}
\end{equation}
In other words, for every $\delta>0$ there is an $N$, such that for every $n>N$
\begin{equation*}\label{eq:new.var.11}
\begin{aligned}
    P\left(\min_{I\subset [n], |I|=k}\vert \M_I\vert >\frac{k!}{n^{k-\varepsilon}}\cdot\dfrac{1}{\sqrt{k}}\right)
    &< \delta.
\end{aligned}
\end{equation*}
Our claim follows from this inequality.

\qedForProof
\end{Proof}

\begin{Proof}{\ of Theorem~\ref{thm:general-result}}
By assumption, there exist $\alpha>0$ and $\beta>0$ such that $f(x)\ge \beta$ for all $x\in[\mu-\alpha,\mu+\alpha]$.
Replacing the variables $X_i$ by $(X_i-\mu)/\alpha$,
we may assume that $\mu=0$ and $\alpha=1$.
We still cannot use Lemma~\ref{lem:Uniform-Result}
since the sample is not from $U(-1,1)$.

To obtain a sample from $U(-1,1)$,
we remove observations as follows:
\begin{itemize}
\item Each $X_i$ outside $[-1,1]$ is removed.
\item If $X_i\in[-1,1]$, it is accepted with probability $\beta/f(X_i)$.
\end{itemize}

Clearly, the surviving sub-sample is $U(-1,1)$-distributed.
(Interestingly, this is a kind of von~Neumann rejection sampling algorithm in reverse.)

Let $\delta,\varepsilon>0$.
We apply Lemma~\ref{lem:Uniform-Result}
with $\delta'=\delta/2$ and $\varepsilon'=\varepsilon$.
Let $N'$ be the constant guaranteed by that lemma, and define
\[
N=\max\left\{\frac{N'}{\beta},\frac{4\log(2/\delta)}{\beta}\right\}.
\]

Let $n\ge N$.
The probability that $X_i$ is accepted is $2\beta$,
and hence $Y\sim B(n,2\beta)$ denotes
the number of accepted observations.
By Chernoff’s bound \cite[Thm.~4.5]{mitzenmacher2017probability},
\begin{equation*}
P(Y\le n\beta)\le e^{-n\beta/4}.
\end{equation*}

If $Y>N'$, then by Lemma~\ref{lem:Uniform-Result},
with probability at least $1-\delta/2$
there exists a $k$-subset whose mean lies in
\[
\left[
-\frac{k!}{n^{k-\varepsilon}\sqrt{k}},
\frac{k!}{n^{k-\varepsilon}\sqrt{k}}
\right].
\]
Hence,
\[
P\!\left(
\exists\, I\subseteq[n],\ |I|=k :
|M_I|
\le
\frac{k!}{n^{k-\varepsilon}\sqrt{k}}
\right)
\ge(1-\delta/2)^2\ge1-\delta .
\]
\qedForProof
\end{Proof}

\begin{remark}
The removal of observations is not meant to describe an algorithmic restriction.
The teacher may benefit from using extreme points;
we only show that the remaining observations suffice.
\end{remark}

\begin{Proof}{\ of Proposition~\ref{prop:scale-general}}
Applying Theorem~\ref{thm:general-result} to the 
transformed 
sample $(Y_i)_{i=1}^n$ yields that for every $\delta>0$ and $\varepsilon>0$, and for all sufficiently large $n$,
\begin{equation}
\label{eq:smooth-proof-mean}
P\!\left(
\exists\, I\subseteq[n],\ |I|=k :
\left|
\frac{1}{k}\sum_{i\in I} Y_i - m
\right|
\le
\frac{k!}{n^{k-\varepsilon}\sqrt{k}}
\right)
\ge
1-\delta .
\end{equation}
Let $I\subseteq[n]$ be a subset whose existence is guaranteed (with probability $1-\delta$) by
\eqref{eq:smooth-proof-mean}.

Take a small interval $J$ around $m$, in which $g'$ is continuous. Let $C$ be an upper bound on $|g'|$ in this interval.
By the mean‑value theorem, there exists a point $\xi_I$
lying between $\frac{1}{k}\sum_{i\in I} Y_i$ and $m$ such that
\begin{equation}
\label{eq:smooth-proof-mvt}
\hat{\theta}_I-\theta^*
=
g\!\left(
\frac{1}{k}\sum_{i\in I} Y_i
\right)-g(m)
=
g'(\xi_I)
\left(
\frac{1}{k}\sum_{i\in I} Y_i - m
\right).
\end{equation}
If $n$ is sufficiently large the average
$\frac{1}{k}\sum_{i\in I} Y_i$ lies in $J$,
and therefore $|g'(\xi_I)|$ is bounded by~$C$.
Combining this bound with \eqref{eq:smooth-proof-mean}
and \eqref{eq:smooth-proof-mvt} yields
\[
\left|
\hat{\theta}_I-\theta^*
\right|
\le
C \cdot \frac{k!}{n^{k-\varepsilon}\sqrt{k}}
\]
with probability at least $1-\delta$. 
Since the proposition is stated for arbitrary $\varepsilon>0$,
the constant $C$ may be absorbed into the rate by replacing
$\varepsilon$ with a smaller positive value.
\qedForProof
\end{Proof}

\begin{Proof}{ of Theorem \ref{thm:median}}
Subtracting the median from every observation, we may assume that
$M=0$.

Apply Theorem~\ref{thm:general-result} with $k=2$.
Since the density is bounded away from $0$ in a neighborhood of $0$,
the assumptions of that theorem are satisfied. Therefore, for every
$\delta>0$ and $\varepsilon>0$, with probability at least $1-\delta$
and for all sufficiently large $n$, there exist distinct indices
$i,j\in[n]$ such that
\[
\left|
\frac{X_i+X_j}{2}
\right|
\le
\frac{2!}{n^{2-\varepsilon}\sqrt2}.
\]

Since the median of the two-point sample $\{X_i,X_j\}$ is
$(X_i+X_j)/2$, the learner can be taught this value by presenting the
subset $\{X_i,X_j\}$.
Absorbing the constant $2!/\sqrt2$ into the rate completes the proof.

\qedForProof
\end{Proof}

\begin{Proof}{\ of Theorem~\ref{thm:median}}
Similarly to the proof of Theorem~\ref{thm:general-result},
let $\alpha,\beta>0$ be such that
$f(x)\ge\beta$ for all $x\in[M-\alpha,M+\alpha]$.
Retaining each observation in this interval independently with
probability $\beta/f(X_i)$ and discarding all remaining observations,
we obtain an i.i.d.\ sample from the uniform distribution on
$[M-\alpha,M+\alpha]$.

If $Y$ denotes the number of surviving observations, then, as in the
proof of Theorem~\ref{thm:general-result},
$Y\sim B(n,2\alpha\beta)$, and therefore $Y=\Omega(n)$ w.h.p.

For the uniform distribution on $[M-\alpha,M+\alpha]$, both the mean
and the median are equal to $M$. Applying
Theorem~\ref{thm:general-result} with $k=2$, we obtain that for every
$\delta>0$ and $\varepsilon>0$, with probability at least $1-\delta$
and for all sufficiently large $n$, there exist two surviving
observations $X_i$ and $X_j$ such that
\[
\left|
\frac{X_i+X_j}{2}-M
\right|
\le
\frac{\sqrt 2}{Y^{\,2-\varepsilon}}.
\]

Since $Y=\Omega(n)$ with high probability, it follows that w.h.p.
\[
\left|
\frac{X_i+X_j}{2}-M
\right|
=
O\!\left(n^{-2+\varepsilon}\right).
\]
Since the median of the two-point sample
$\{X_i,X_j\}$ equals $(X_i+X_j)/2$, the result follows.
\qedForProof
\end{Proof}

\bibliography{references/references}

\end{document}